\documentclass{article}

\usepackage[heads=vee,repositionpullbacks,midshaft,PostScript=dvips,nohug,scriptlabels,size=2em]{diagrams}
\usepackage{bbm}
\usepackage{amssymb,amsmath}
\usepackage{mathrsfs}
\usepackage{color}

\usepackage{natbib,natbibspacing}
\bibpunct{(}{)}{,}{a}{}{,}
\renewcommand{\bibfont}{\small}

\definecolor{myurlcolor}{rgb}{0,0,0.5}
\usepackage{hyperref}
\hypersetup{colorlinks,
linkcolor=black,
citecolor=black,
urlcolor=myurlcolor}

\newcommand{\emptybk}{\:\:}
\newcommand{\blank}{(\emptybk)}
\newcommand{\dashbk}{-}
\newcommand{\cat}[1]{\mathscr{#1}}
\newcommand{\fcat}[1]{\mathbf{#1}}
\newcommand{\such}{\:|\:}
\newcommand{\Hom}{\mathrm{Hom}}
\newcommand{\op}{\mathrm{op}}
\newcommand{\Set}{\fcat{Set}}
\newcommand{\goesto}{\mapsto}
\newcommand{\parpair}[2]{\pile{\rTo^{\scriptstyle #1}\\ 
\rTo_{\scriptstyle #2}}}
\newcommand{\parpairu}{\rightrightarrows}
\newcommand{\vslob}[3]
	{\left.
	\begin{diagram}[height=1.5em]
	#1		\\
	\dTo>{\,#2}	\\
	#3		\\
	\end{diagram}
	\right.}
\newarrow{Incl}C--->
\newarrow{Mono}>--->
\newcommand{\reals}{\mathbb{R}}
\newlength{\templength}    
\newcommand{\demph}[1]{\textbf{\textup{#1}}}
\newcommand{\scat}[1]{\mathbb{#1}}
\newcommand{\iso}{\cong}
\newcommand{\nat}{\mathbb{N}}	
\newcommand{\eqv}{\simeq}
\newcommand{\pr}{\mathrm{pr}}
\newcommand{\of}{\,\raisebox{0.08ex}{\ensuremath{\scriptstyle\circ}}\,}
\newcommand{\sub}{\subseteq}
\newcommand{\cell}[4]{\put(#1,#2){\makebox(0,0)[#3]{\ensuremath{#4}}}}
\newcommand{\Colt}[1]{{\displaystyle\lim_{\rightarrow #1}}\,}
\newcommand{\tr}{\texttt{t}}
\newcommand{\fa}{\texttt{f}}
\newcommand{\One}{\mathbbm{1}}
\newcommand{\chunk}[1]{\subsection*{#1}}
\newcommand{\slogan}[1]{\begin{center}\it #1\end{center}}
\newcommand{\E}{\cat{E}}
\newcommand{\F}{\cat{F}}
\newcommand{\T}{\cat{T}}
\newcommand{\Mono}{\fcat{Mono}}
\newcommand{\Sub}{\mathrm{Sub}}
\newcommand{\toby}[1]{\stackrel{#1}{\to}}
\newcommand{\Pshf}[1]{\fcat{Psh}{#1}}
\newcommand{\Psh}[1]{\widehat{#1}}
\newcommand{\Sh}{\fcat{Sh}}
\newcommand{\Et}[1]{\fcat{Et}{#1}}
\newcommand{\Open}{\fcat{Open}}
\newcommand{\FinSet}{\fcat{FinSet}}
\newcommand{\restr}[1]{\vert_{#1}}
\newcommand{\TopSp}{\fcat{TopSp}}
\newcommand{\Topos}{\fcat{Topos}}
\newcommand{\Frame}{\fcat{Frame}}
\newcommand{\Loc}{\fcat{Loc}}
\newcommand{\incl}{\hookrightarrow}
\newcommand{\fp}{\mathrm{fp}}
\newcommand{\Gp}{\fcat{Grp}}
\newcommand{\monic}{\rightarrowtail}
\newcommand{\monicby}[1]{\stackrel{#1}{\monic}}
\newcommand{\epic}{\twoheadrightarrow}
\newcommand{\dbot}{{\scriptstyle\bot}}
\newcommand{\deqv}{{\scriptstyle\eqv}}
\newcommand{\cln}{\colon}

\newtheorem{thm}{Theorem}[section]
\newtheorem{lemma}[thm]{Lemma}
\newtheorem{fact}[thm]{Fact}
\newtheorem{predefn}[thm]{Definition}
\newenvironment{defn}{\begin{predefn}\upshape}{\end{predefn}}
\newtheorem{preexamples}[thm]{Examples}
\newenvironment{examples}{\begin{preexamples}\upshape}{\end{preexamples}}

\title{An informal introduction to topos theory} 
\author{Tom Leinster%
\thanks{School of Mathematics and Statistics, University of Glasgow, Glasgow
G12 8QW, UK;
Tom.Leinster@glasgow.ac.uk.  Supported by an EPSRC Advanced Research
Fellowship.}}
\date{}

\begin{document}

\maketitle

\vspace*{10mm}
\tableofcontents
\vspace*{12mm}

\noindent
This short text is for readers who are confident in basic category theory but
know little or nothing about toposes.  It is based on some impromptu talks
given to a small group of category theorists.  I am no expert on topos theory.
These notes are for people even less expert than me.

In keeping with the spirit of the talks, what follows is light on both detail
and references.  For the reader wishing for more, almost everything
here is presented in respectable form in Mac Lane and Moerdijk's very pleasant
introduction to topos theory (\citeyear{MaMo}).  Nothing here is new, not even
the expository viewpoint (very loosely inspired by \citet{JohSE}).

As a rough indication of the level of knowledge assumed, I will take it that
you are totally comfortable with the Yoneda Lemma and the concept of
cartesian closed category, but I will not assume that you know the definition
of subobject classifier or of topos.

Section~\ref{sec:defn} explains the definition of topos.  The remaining three
sections discuss some of the connections between topos theory and other
subjects.  There are many more such connections than I will mention; I hope it
is abundantly clear that these notes are, by design, a quick sketch of a
large subject.

Section~\ref{sec:sets} is on connections between topos theory and set theory.
There are two themes here.  One is that, using the language of toposes, we can
write down an axiomatization of sets that sticks closely to how sets are
actually used in mathematics.  This provides an appealing alternative to ZFC.
The other, related, theme is that 
\slogan{a topos is a generalized category of sets.}

Section~\ref{sec:geom} is on connections with geometry (in a broad sense);
there the thought is that
\slogan{a topos is a generalized space.}

Section~\ref{sec:univ-alg} is on connections with universal algebra:
\slogan{a topos is a generalized theory.}
What this means is that there is one topos embodying the concept of `ring',
another embodying the concept of `field', and so on.  This is the story of
classifying toposes.

\setlength{\templength}{\parindent}
\noindent
\parbox{.7\textwidth}{%
\setlength{\parindent}{\templength}
Sections~\ref{sec:sets}--\ref{sec:univ-alg} can be read in any order, except
that ideally~\S\ref{sec:geom} (geometry) should come
before~\S\ref{sec:univ-alg} (universal algebra).  You \emph{can}
read~\S\ref{sec:univ-alg} without having read~\S\ref{sec:geom}, but the price
to pay is that the notion of `geometric morphism'---defined in~\S\ref{sec:geom}
and used in~\S\ref{sec:univ-alg}---might seem rather
mysterious.}%
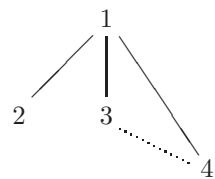
\parbox{.05\textwidth}{\ \ \ \ }
\parbox{.25\textwidth}{%
\setlength{\unitlength}{1em}
\begin{picture}(8,7.5)(1.2,1)
\cell{5}{8}{t}{1}
\cell{1.7}{4}{c}{2}
\cell{5}{4}{c}{3}
\cell{8.8}{2}{c}{4}
\put(4.5,7){\line(-1,-1){2.3}}
\put(5,7){\line(0,-1){2.3}}
\put(5.5,7){\line(2,-3){3}}
\qbezier[12](5.5,3.5)(6.8,2.85)(8.1,2.2)
\end{picture}
}

Algebraic geometers beware: the word `topos' is used by mathematicians in two
slightly different senses, according to circumstance and culture.  There are
elementary toposes and Grothendieck toposes.  Category theorists tend to use
`topos' to mean `elementary topos' by default, although Grothendieck toposes
are also important in category theory.  But when an algebraic geometer says
`topos', they almost certainly mean `Grothendieck topos' (what else?).

Grothendieck toposes are categories of sheaves.  Elementary toposes are
slightly more general, and the definition is simpler.  They are what I will 
emphasize here.  Grothendieck toposes are the subject of
Section~\ref{sec:geom}, and appear fleetingly elsewhere; but if you only want
to learn about categories of sheaves, this is probably not the text for you.

\paragraph*{Acknowledgements}  I thank Andrei Akhvlediani, Eugenia Cheng,
Richard Garner, Nick Gurski, Ignacio Lopez Franco and Emily Riehl for their
participation and encouragement.  Aspects of Section~\ref{sec:univ-alg} draw
on a vaguely similar presentation of vaguely similar material by Richard
Garner.  I thank the organizers of Category Theory 2010 for making the talks
possible, even though they did not mean to: Francesca Cagliari, Eugenio Moggi,
Marco Grandis, Sandra Mantovani, Pino Rosolini, and Bob Walters.  I thank Jon
Phillips, Urs Schreiber, Mike Shulman, Alex Simpson, Danny Stevenson and Todd
Wilson for suggestions and corrections.  I am especially grateful to Todd
Trimble for carefully reading an earlier version and suggesting many
improvements.  The commutative diagrams were made using Paul Taylor's macros.

\section{The definition of topos}
\label{sec:defn}

The hardest part of the definition of topos is the concept of subobject
classifier, so I will begin there.  For motivation, I will speak of `the
category of sets' (and functions).  What exactly this means will be discussed
in Section~\ref{sec:sets}, but for now we proceed informally.

In the category of sets, inverse images are a special case of pullbacks.
That is, given a map $f\cln X \to Y$ of sets and a subset $B \sub Y$, we have a
pullback square
\[
\begin{diagram}
f^{-1}B\SEpbk   &\rTo   &B      \\
\dIncl          &       &\dIncl \\
X               &\rTo_f &Y.     \\
\end{diagram}
\]
In particular, this holds when $B$ is a 1-element subset $\{y\}$ of $Y$:
\[
\begin{diagram}
f^{-1}\{y\}\SEpbk         &\rTo   &\{y\}  \\
\dIncl                  &       &\dIncl \\
X                       &\rTo_f &Y.     \\
\end{diagram}
\]
There is no virtue in distinguishing between one-element sets, so we might as
well write $1$ instead of $\{y\}$; then the inclusion $\{y\} \incl Y$ becomes
the map $1 \to Y$ picking out $y \in Y$, and we have a pullback square
\[
\begin{diagram}
f^{-1}\{y\}\SEpbk       &\rTo^! &1      \\
\dIncl                  &       &\dTo>y \\
X                       &\rTo_f &Y.     \\
\end{diagram}
\]

Next consider characteristic functions of subsets.  Fix a two-element set $2 =
\{\tr, \fa\}$ (`true' and `false').  Then for any set $X$, the subsets of $X$
are in bijective correspondence with the functions $X \to 2$.  In one
direction, given a subset $A \sub X$, the corresponding function $\chi_A\cln X
\to 2$ is defined by
\[
\chi_A(x)
=
\begin{cases}
\tr     &\text{if } x \in A     \\
\fa     &\text{if } x \not\in A
\end{cases}
\]
($x \in X$).  In the other, given a function $\chi\cln X \to 2$, the
corresponding subset of $X$ is $\chi^{-1}\{\tr\}$.  To say that this latter
process $\chi \mapsto \chi^{-1}\{\tr\}$ is a bijection is to say that for all
$A \sub X$, there is a unique function $\chi\cln X \to 2$ such that $A =
\chi^{-1}\{\tr\}$.  In other words: for all $A \sub X$, there is a unique
function $\chi\cln X \to 2$ such that
\[
\begin{diagram}
A       &\rTo^!         &1              \\
\dIncl  &               &\dTo>\tr       \\
X       &\rTo_\chi      &2              \\
\end{diagram}
\]
is a pullback square.

This property of sets can now be stated in purely categorical terms.  We use
$\monic$ to indicate a mono ($=$ monomorphism $=$ monic).  

\begin{defn}
Let $\E$ be a category with a terminal object, $1$.  A \demph{subobject
classifier} in $\E$ is an object $\Omega$ together with a map $\tr\cln 1 \to
\Omega$ such that for every mono $A \monicby{m} X$ in $\E$, there exists a
unique map $\chi\cln X \to \Omega$ such that
\[
\begin{diagram}
A       &\rTo^!         &1              \\
\dMono<m&               &\dTo>\tr       \\
X       &\rTo_\chi      &\Omega         \\
\end{diagram}
\]
is a pullback square.
\end{defn}

So, we have just observed that $\Set$ has a subobject classifier, namely, the
two-element set.  In the general setting, we may write $\chi$ as $\chi_A$ (or
properly, $\chi_m$) and call it the \demph{characteristic function} of $A$ (or
$m$).

To understand this further, we need two lemmas.

\begin{lemma}   \label{lemma:pb-mono}
In any category, the pullback of a mono is a mono.  That is, if
\[
\begin{diagram}
\cdot           &\rTo   &\cdot  \\
\dTo<{m'}       &       &\dTo>m \\
\cdot           &\rTo   &\cdot  \\
\end{diagram}
\]
is a pullback square and $m$ is a mono, then so is $m'$.  
\end{lemma}

\begin{lemma}
In any category with a terminal object $1$, every map out of $1$ is a mono.
\end{lemma}

So, pulling $\tr\cln 1 \to \Omega$ back along \emph{any} map $X \to \Omega$
gives a mono into $X$.

It will also help to know the result of the following little exercise.  It
says, roughly, that in the definition of subobject classifier, the fact that
$1$ is terminal comes for free.

\begin{fact}    \label{fact:one-terminal}
Let $\E$ be a category and let $T \monicby{\tr} \Omega$ be a mono in $\E$.
Suppose that for every mono $A \monicby{m} X$ in $\E$, there is a unique map
$\chi\cln X \to \Omega$ such that there is a pullback square
\[
\begin{diagram}
A               &\rTo           &T              \\
\dMono<m        &               &\dMono>\tr     \\
X               &\rTo_\chi      &\Omega.        \\
\end{diagram}
\]
Then $T$ is terminal in $\E$.
\end{fact}

This leads to a second description of subobject classifiers.  Let $\Mono(\E)$
be the category whose objects are monos in $\E$ and whose maps are pullback
squares.  Then a subobject classifier is exactly a terminal object of
$\Mono(\E)$.  

Here is a third way of looking at subobject classifiers.  Given a category
$\E$ and an object $X$, a \demph{subobject} of $X$ is officially an
isomorphism class of monos $A \monicby{m} X$ (where isomorphism is taken in the
slice category $\E/X$).  For example, when $\E = \Set$, two monos
\[
A \monicby{m} X, 
\quad
A' \monicby{m'} X
\]
are isomorphic if and only if they have the same image; so subobjects of $X$
correspond one-to-one with subsets of $X$.  I say `officially' because half
the time people use `subobject of $X$' to mean simply `mono into $X$', or slip
between the two meanings without warning.  It is a harmless abuse of
language, which I will adopt.   

For $X \in \E$, let $\Sub(X)$ be the class of subobjects (in the official
sense) of $X$.  Assume that $\E$ is well-powered, that is, each $\Sub(X)$ is a
set rather than a proper class.  Assume also that $\E$ has pullbacks.  By
Lemma~\ref{lemma:pb-mono}, every map $X \toby{f} Y$ in $\E$ induces a map
$\Sub(Y) \toby{f^*} \Sub(X)$ of sets, by pullback.  This defines a functor
$\Sub\cln \E^\op \to \Set$.

Third description: a subobject classifier is a representation of this functor
$\Sub$.

This makes intuitive sense, since for $\Sub$ to be representable means that
there is an object $\Omega \in \E$ satisfying
\[
\Sub(X) \iso \E(X, \Omega)
\]
naturally in $X \in \E$.  In the motivating case of the category of sets, this
directly captures the thought that subsets of a set $X$ correspond naturally
to maps $X \to \{\tr, \fa\}$.  

Now we show that this is equivalent to the original definition.  By the Yoneda
Lemma, a representation of $\Sub\cln \E^\op \to \Set$ amounts to an object
$\Omega \in \E$ together with an element $\tr \in \Sub(\Omega)$ that is
`generic' in the following sense:
\begin{quote}
for every object $X \in \E$ and element $m \in \Sub(X)$, there is a unique map
$\chi\cln X \to \Omega$ such that $\chi^*(\tr) = m$.
\end{quote}
In other words, a representation of $\Sub$ is a mono $T \monicby{\tr} \Omega$
in $\E$ satisfying the condition in Fact~\ref{fact:one-terminal}.  In other
words, it is a subobject classifier.

\begin{defn}
A \demph{topos} (or \demph{elementary topos}) is a cartesian closed category
with finite limits and a subobject classifier.
\end{defn}

\begin{examples}        \label{egs:toposes}
\begin{enumerate}
\item  
The primordial topos is $\Set$.  It has special properties not shared by most
other toposes.  This is the subject of Section~\ref{sec:sets}.

\item 
For any set $I$, the category $\Set^I$ of $I$-indexed families of sets is a
topos.  Its subobject classifier is the constant family $(2)_{i \in I}$, where
$2$ is a two-element set.

\item 
For any group $G$, the category $\Set^G$ of left $G$-sets is a topos.  Its
subobject classifier is the set $2$ with trivial $G$-action.

\item   \label{eg:topos-pshf}
Encompassing all the previous examples, if $\scat{A}$ is any small category
then the category $\Psh{\scat{A}} = \Set^{\scat{A}^\op}$ of presheaves on
$\scat{A}$ is a topos.  We can discover what its subobject classifier must be
by a thought experiment: \emph{if} $\Omega$ is a subobject classifier then by
the Yoneda Lemma,
\[
\Omega(a) 
\iso
\Psh{\scat{A}}( \scat{A}(\dashbk, a), \Omega )
\iso
\Sub(\scat{A}(\dashbk, a))
\]
for all $a \in \scat{A}$.  So $\Omega(a)$ must be the set of subfunctors of
$\scat{A}(\dashbk, a)$; and one can check that defining $\Omega(a)$ in this
way does indeed give a subobject classifier.  A subfunctor of
$\scat{A}(\dashbk, a)$ is called a \demph{sieve} on $a$; it is a collection of
maps into $a$ satisfying a certain condition. 

\item For any topological space $S$, the category $\Sh(S)$ of sheaves on $S$
is a topos.  This is the subject of Section~\ref{sec:geom}.  Modulo a small
lie that I will come back to there, the space $S$ can be recovered from the
topos $\Sh(S)$.  Hence the class of spaces embeds into the class of toposes,
and this is why toposes can be viewed as generalized spaces.

Sheaves will be defined and explained in Section~\ref{sec:geom}.  To give a
brief sketch: denote by $\Open(S)$ the poset of open subsets of $S$; then a
\demph{presheaf} on the space $S$ is a presheaf on the category $\Open(S)$,
and a sheaf on $S$ is a presheaf with a further property.  I will
consistently use `sheaf' to mean what some would call `sheaf of sets'.  A
sheaf of groups, rings, etc.\ is the same as an internal group, ring etc.\ in
$\Sh(S)$.  

\item
The category $\FinSet$ of finite sets is a topos.  Similarly, $\Set$ can be
replaced by $\FinSet$ in all of the previous examples, giving toposes of
finite $G$-sets, finite sheaves, etc.  
\end{enumerate}
\end{examples}

You might ask `why is the definition of topos what it is?  Why that
\emph{particular} collection of axioms?  What's the motivation?'  I will not
attempt to answer, except by explaining several ways in which the definition
has been found useful.  It is also worth noting that the topos axioms have
many non-obvious consequences, giving toposes a far richer structure than most
categories.  For example, every map in a topos factorizes, essentially
uniquely, as an epi followed by a mono.  More spectacularly, the axioms imply
that every topos has finite \emph{co}limits.  This can be proved by the
following very elegant strategy, due to \citet{Pare}.  For every topos $\E$,
we have the contravariant power set functor $P = \Omega^{(\dashbk)}\cln \E^\op
\to \E$.  It can be shown that $P$ is monadic.  But monadic functors create
limits, and $\E$ has finite limits.  Hence $\E^\op$ has finite limits; that
is, $\E$ has finite colimits.

\section{Toposes and set theory}
\label{sec:sets}

Here I will describe what makes `the' category of sets special among all
toposes, and explain why I just put `the' in quotation marks.  This is the
stuff of revolution: it can completely change your view of set theory.  It
also provides an invaluable insight into topos theory as a whole.

We begin by listing some special properties of the topos $\Set$, using only
the most commonplace assumptions about how sets and functions behave.

\begin{enumerate}
\item[\textbf{1.}]  
The terminal object $1$ is a separator (generator).  That is, given maps $X
\parpair{f}{g} Y$ in $\Set$, if $f \of x = g \of x$ for all $x\cln 1 \to X$
then $f = g$.

It is worth dwelling on what this says.  Maps $1 \to X$ correspond to elements
of $X$, and we make no notational distinction between the two.  Moreover,
given an element $x \in X$ and a map $f\cln X \to Y$, we can compose the maps
\[
1 \toby{x} X \toby{f} Y
\]
to obtain a map $f \of x\cln 1 \to Y$, and this is the map corresponding to the
element $f(x) \in Y$.  (We might harmlessly write both $f \of x$ and $f(x)$ as
$fx$.)  Thus, elements are a special case of functions, and evaluation is a
special case of composition.

The property above says that if $f(x) = g(x)$ for all $x \in X$ then $f = g$.
In other words, a function is determined by its effect on elements.

\item[\textbf{2.}]
Write $0$ for the initial object of $\Set$ (the empty set).  Then $0 \not\iso
1$.  Equivalently, $\Set$ is not equivalent to the terminal category $\One$.
\end{enumerate}

A topos satisfying properties~\textbf{1} and~\textbf{2} is called
\demph{well-pointed}.  

\begin{enumerate}
\item[\textbf{3.}] 
This property says, informally, that there is a set consisting of the natural
numbers. 

What are the `the natural numbers', though?  One way to get at an answer is
to use the principle that sequences can be defined recursively.  That is,
given a set $X$, an element $x \in X$, and a map $r\cln X \to X$, there is a
unique sequence $(x_n)_{n = 0}^\infty$ in $X$ such that 
\begin{equation}
\label{eq:recursion}
x_0 = x,
\quad
x_{n + 1} = r(x_n) 
\quad(n \in \nat).
\end{equation}
A sequence $(x_n)_{n = 0}^\infty$ in $X$ is just a map $f\cln \nat \to X$,
and if we write $s\cln \nat \to \nat$ for the function $n \mapsto n + 1$
(`successor'), then~\eqref{eq:recursion} says exactly that the diagram
\begin{equation}
\label{eq:nat}
\begin{diagram}
        &               &\nat   &\rTo^s &\nat   \\
1       &\ruTo(2,1)<0   &\dTo>f &       &\dTo>f \\
        &\rdTo(2,1)<x   &X      &\rTo_r &X      \\
\end{diagram}
\end{equation}
commutes.
\end{enumerate}

\begin{defn}
Let $\E$ be a category with a terminal object, $1$.  A \demph{natural numbers
object} in $\E$ is a triple $(N, 0, s)$, with $N \in \E$, $0\cln 1 \to N$, and
$s\cln N \to N$, that is initial as such: for any triple $(X, x, r)$ of
the same type, there is a unique map $f\cln N \to X$ such that~\eqref{eq:nat}
commutes (with $N$ in place of $\nat$).
\end{defn}

Property~\textbf{3} is, then, that $\Set$ has a natural numbers object.

\begin{enumerate}
\item[\textbf{4.}] 
Epis split.  That is, for any epimorphism (surjection) $e\cln X \to Y$ in
$\Set$, there exists a map $m\cln Y \to X$ such that $e \of m = 1_Y$.  The
splitting $m$ chooses for each $y \in Y$ an element of the nonempty set
$e^{-1}\{y\}$.  The existence of such splittings is precisely the Axiom of
Choice.  Generally, a category is said to satisfy the \demph{Axiom of Choice}
(or to `have Choice') if epis split.
\end{enumerate}

In summary,
\slogan{sets and functions form a well-pointed topos\\
with natural numbers object and Choice.}
The category of sets has many other elementary properties (such as the fact
that the subobject classifier has exactly two elements), but they are all
consequences of the properties just mentioned.

But what is this thing called `the category of sets'?  What do we have to
assume about sets in order to prove that these properties hold?  

Many mathematicians do not like to be bothered with such questions, because
they know that the standard answer will be something like `sets are anything
satisfying the axioms of ZFC'---and they feel that ZFC is irrelevant to what
they do, and prefer not to hear about it.

The standard answer is \emph{valid}, in the sense that for every model of
ZFC, there is a resulting category of sets satisfying the properties
above.  But it may seem \emph{irrelevant}, because at no point in establishing
the properties did it feel necessary to call on an axiom system: all the
properties are suggested directly by the naive imagery of a set as a bag of
dots.  

There is, however, another type of answer---and this was Lawvere's radical
idea.  It is this: 
\slogan{we take the properties above as our axioms on sets.}
In other words, we do away with ZFC entirely, and ask instead that sets and
functions form a well-pointed topos with natural numbers object and Choice.
`The' category of sets is any category satisfying these axioms.  In fact we
should say \emph{a} category of sets, since there may be many different such
categories, as we shall see.

This is Lawvere's Elementary Theory of the Category of Sets (ETCS), stated in
modern language.  (See \citet{LawETCS}, or \citet{LaRo} for a good expository
account.)  It is nearly fifty years old, but still has not gained the currency
it deserves, for reasons on which one can speculate.

\paragraph*{Digression} You might be thinking that this is circular: that this
axiomatization of sets 
depends on the notion of category, and the notion of category depends on some
notion of collection or set.  But in fact, ETCS does not depend on the general
notion of category.  It can be stated without using the word `category' once.

To see this, we need to back up a bit.  The ZFC axiomatization of sets looks,
informally, like this:
\begin{itemize}
\item there are some things called `sets'
\item there is a binary relation `$\in$' on sets
\item some axioms hold.
\end{itemize}
People seeing this (or the formal version) often ask certain questions.  What
does `some things' mean?  Do you mean that there is a \emph{set} of sets?
(No.)  What exactly is meant by `binary relation'?  (It means that for each
set $X$ and set $Y$, the statement `$X \in Y$' is deemed to be either true or
false.)  What do you mean, `deemed'?  Etc.  This is not a logic course, and I
will not attempt to answer the questions except to say that there is an
assumed common understanding of these terms.  To hide behind jargon, ZFC is a
first-order theory.

The ETCS axiomatization of sets looks like this:
\begin{itemize}
\item there are some things called `sets'
\item for each set $X$ and set $Y$, there are some things called `functions
from $X$ to $Y$'
\item for each set $X$, set $Y$ and set $Z$, there is a binary operation
assigning to each pair of functions
\[
f\cln X \to Y, 
\quad
g\cln Y \to Z
\]
a function $g \of f\cln X \to Z$
\item some axioms hold.
\end{itemize}
You can ask the same kind of logical questions as for ZFC---what exactly is
meant by `binary operation'?\ etc.---which again I will not attempt to
answer.  The difficulties are no worse than for ZFC, and again, in the jargon,
ETCS is a first-order theory.  

Stated in this way, the ETCS axioms begin by saying that composition is
associative and has identities (so that sets, functions and composition of
functions define a category); then they say that binary products and
equalizers of sets exist, and there is a terminal set (so that the category of
sets has finite limits); and so on, until we have said that sets and functions
form a well-pointed topos with natural numbers object and Choice.  You can do
it in about ten axioms.

Here ends the digression.

\bigskip

ZFC axiomatizes sets and membership, whereas ETCS axiomatizes sets and
functions.  Anything that can be expressed in one language can be expressed in
the other: in the usual implementation of ZFC, a function $X \to Y$ is
defined as a suitable subset of $X \times Y$, and in ETCS, an element of $X$
is defined as a function from the terminal set to $X$.  But an advantage of
the categorical approach is that it avoids the chains of elements of elements
of elements that are so important in traditional set theory, yet
seem so distant from most of mathematics.

ZFC is slightly stronger than ETCS.  `Stronger' means that everything that can
be deduced about sets from the ETCS axioms can also be deduced in ZFC, but not
vice versa.  `Slightly' is meant in a sociological sense.  I believe it has
been said that the mathematics in an ordinary undergraduate syllabus
(excluding, naturally, any course in ZFC) makes no more assumptions about sets
than are made by ETCS.  If that is so, it must also be the case that for many
mathematicians, nothing in their entire research career requires more than
ETCS.

The technical relationship between ZFC and ETCS is well understood.  It is
known exactly which fragment of ZFC is equivalent to ETCS (namely, `bounded'
or `restricted' Zermelo with Choice; see \citet{MaMo}).  It is also known what
needs to be added to ETCS in order to obtain a system of equal strength to
ZFC.  This extra ingredient is an axiom scheme (a countably infinite family of
axioms) that set theorists in the traditional mould would call Replacement,
and category theorists would call a form of cocompleteness.  It says,
informally, that given any set $I$ and family $(X_i)_{i \in I}$ of sets
specified by a first-order formula, the coproduct $\sum_{i \in I} X_i$ exists.
The existence of this coproduct is expressed by saying that there exist a set
$X$ and a map $p\cln X \to I$ (to be thought of as the projection $\sum_{i \in
I} X_i \to I$) such that for each $i \in I$, the inverse image $p^{-1}\{i\}$
is isomorphic to $X_i$.  See Section~8 of \citet{McLECS} for details.

Topos theory therefore provides a different viewpoint on set theory.  Let us
take a brief look from this new viewpoint at a famous theorem of set theory:
that the Continuum Hypothesis is independent of the usual set-theoretic
axioms, as proved by G\"odel and Cohen.

Temporarily, let us say that a `category of sets' is a well-pointed topos with
natural numbers object and Choice, satisfying the axiom scheme of Replacement.
A category of sets is said to \demph{satisfy the Continuum Hypothesis} if for
all objects $X$,
\begin{align*}
                &\text{there exist monos } N \monic X \monic 2^N        \\
\implies        &X \iso N \text{ or } X \iso 2^N.
\end{align*}
(As usual, $N$ denotes the natural numbers object; $2$ is the subobject
classifier.)  Stated categorically, the theorem is this: given any category of
sets, you can build one that satisfies the Continuum Hypothesis and one that
does not.  This is only a rephrasing of the standard statement, but if you are
more at home with the term `category' than with `model of a first-order
theory', you might find it less mysterious.

So far we have seen the benefits of viewing the/a category of sets as a
special topos.  But the other way round, there are great benefits to viewing a
topos as a generalized category of sets.  For example, we might view
$\Set^\nat$ as the category of sets varying through (discrete) time.  The set
of human beings alive today is an object of $\Set^\nat$: as the meaning of
`today' changes, the set changes.  A sheaf can similarly be understood as a
set varying through space.

People (especially Lawvere) sometimes refer to the category of sets as the (or
a) topos of \emph{constant} sets, to contrast it with toposes of variable
sets.  There are also toposes whose objects can informally be thought of as
`cohesive' sets, which means the following.  In an ordinary set, the points
have no relation or attachment to each other: they do not `cohere'.  But a
cohesive set carries something like a topology or smooth structure, so that
the points are in some sense stuck together.  For example, there are toposes of
smooth spaces, which are the setting for synthetic differential geometry.
From this point of view, the category of ordinary sets is extreme among all
toposes: its objects are sets with no variation or cohesion at all.

Viewing the objects of a topos as generalized sets is much more than a useful
mental technique.  In fact, it is valid to use set-like language and reasoning
in \emph{any} topos, provided that we stick to certain rules.  This language
is called the `internal language' of the topos.  

Many of the central ideas of topos theory are simple, but that simplicity can
easily be obscured by the richness of structure available in a topos.
Such is the case for the internal language.  I will therefore describe the
idea in a much more basic setting.

First let $\E$ be any category whatsoever, and let $A$ be an object of $\E$.
A \demph{generalized element} of $X$ is simply a map in $\E$ with codomain
$X$.  A generalized element $x\cln S \to X$ may be said to be of \demph{shape}
$S$, or to be an \demph{$S$-element} of $X$.  In the special case that $S$ is
terminal, $S$-elements are called \demph{global elements}.  (See
Example~\ref{eg:sheaves}(\ref{eg:shf-extend}) for a hint on the reason for the
name.)  In the category of sets, the global elements are the ordinary
elements, but in other categories, the global elements might be very
uninteresting: consider the category of groups, for instance.

Given a map $f\cln X \to Y$ in $\E$, any generalized element $x$ of $X$ gives
rise to a generalized element $fx$ of $Y$.  This is the composite $f \of x$,
but can also be thought of as `$f(x)$': see the remarks on property~\textbf{1}
at the beginning of this section.  For maps $X \parpair{f}{g} Y$, we have
\[
f = g
\iff
fx = gx \text{ for all generalized elements } x \text{ of } X.
\]
(Proof of $\Leftarrow$: take $x = 1_X$.)  This is emphatically not true if we
replace `generalized' by `global': again, consider groups.  

This language of generalized elements is the \demph{internal language} of the
category.  It fits well with ordinary categorical terminology and notation.
For example, let $\E$ be a category with finite products.  In the internal
language, the definition of product reads, informally: an $S$-element of $X
\times Y$ consists of an $S$-element of $X$ together with an $S$-element of
$Y$.  Apart from the `$S$-' prefixes, this is identical to the ordinary
description of the cartesian product of sets $X$ and $Y$.  And in standard
categorical notation, the map $S \to X \times Y$ with components $x\cln S \to
X$ and $y\cln S \to Y$ is denoted by $(x, y)$, thus extending the set-theoretic
notation for a (global) element of a cartesian product.

To see why the internal language is useful, consider, for instance, internal
groups in a finite product category $\E$.  A group in $\E$ is an object $X$
together with maps
\[
m\cln X \times X \to X,
\quad
i\cln X \to X,
\quad
e\cln 1 \to X
\]
satisfying some axioms.  Those axioms are usually expressed as commutative
diagrams, which have been obtained by translating the classical axioms into
diagrammatic form.  But there is no need to translate them: the classical
axioms can simply be repeated verbatim and interpreted as statements about
\emph{generalized} elements.  This is equivalent.  For example, it is easy to
show that the commutative diagram for associativity is equivalent to the
statement that
\begin{equation}        \label{eq:assoc}
m(m(x, y), z) = m(x, m(y, z))
\end{equation}
for all generalized elements $x, y, z$ of $X$ of the same shape.  (They have
to be the same shape in order for expressions such as $(x, y)$ to make sense.)
And just as for ordinary elements in $\Set$, there is no harm in writing $xy$
instead of $m(x, y)$, and similarly $x^{-1}$ instead of $i(x)$.

More valuably still, \emph{proofs} written down in the classical set-theoretic
scenario will actually be valid in an arbitrary finite product category $\E$,
as long as whatever was said about elements in $\Set$ is also true for
generalized elements in $\E$.  For example, whenever $X$ is a group in $\Set$
and $x, y, a \in X$, we have
\begin{equation}        \label{eq:cancellation}
xa = ya \implies x = y.
\end{equation}
Proof:
\begin{align*}
xa = ya &
\implies
(xa)a^{-1} = (ya)a^{-1}
\implies
x(a a^{-1}) = y(a a^{-1})       \\
&
\implies
xe = ye
\implies
x = y.
\end{align*}
We can immediately conclude that the implication~\eqref{eq:cancellation} holds
whenever $X$ is a group in an arbitrary finite product category $\E$ and $x,
y, a$ are generalized elements of $X$ of the same shape.  Indeed, each step in
the proof is an application of an axiom such as~\eqref{eq:assoc} valid in the
general setting.  

The internal language is a massively labour-saving device.  To prove that an
equation valid in ordinary groups is also valid for internal groups, you
merely need to cast an eye over the proof and convince yourself that it holds
for generalized elements too.  In contrast, try proving the internal version
of the equation
\begin{equation}       \label{eq:inverse}
y^{-1} x^{-1} = (xy)^{-1}
\end{equation}
by diagrammatic methods.  First it has to be \emph{stated} diagrammatically.
It says that the diagram
\[
\begin{diagram}
X \times X      &\rTo^{\mathrm{sym}}    &X \times X     &
\rTo^{i \times i}       &X \times X     \\
\dTo<m          &                       &               &
                        &\dTo>m         \\
X               &                       &\rTo_i         &
                        &X              \\
\end{diagram}
\]
commutes.  Then it has to be \emph{proved}, by filling the
inside of this diagram with instances of the diagrams encoding the group
axioms.  (It seems to need at least ten or so inner diagrams.)  But once you
have an elementwise proof, all this effort is unnecessary.  And the
example~\eqref{eq:inverse} chosen was very simple: for more
complex statements, the benefits of the internal language become clearer
still.

The internal language of toposes is similar to that of finite product
categories, but much richer.  As well as being able to form pairs $(x, y)$ of
generalized elements, we can take generalized elements of exponentials $Y^X$
(to be thought of as families of maps $X \to Y$), form subobjects such as
\[
\{ x \in X \such fx = gx \}
\]
(the equalizer of $X \parpair{f}{g} Y$), and so on.  Almost anything that can
be expressed or proved in the category of sets can be reproduced in an
arbitrary topos.  The only sticking points are the law of the excluded middle
and the axiom of choice.  Any proof that avoids those---any constructive
proof, in a sense that can be made precise---generalizes to an arbitrary
topos.

Phrases with more or less the same meaning as `internal language' are
`Mitchell--B\'enabou language' and `internal logic'.  See, for instance,
\citet{MaMo} or \citet{JohSE}.  There you can also find more spectacular
applications of topos theory to set theory, including topics such as forcing.

\section{Toposes and geometry}
\label{sec:geom}

This section covers concepts such as sheaf, geometric morphism (map of
toposes), Grothendieck topos, and locale.  But the most important thing I
want to explain is how and why geometry has inspired so much of topos theory.

\chunk{Sheaves}

Let $X$ be a topological space.  (Following tradition, I will switch from my
previous convention of using $X$ to denote an object of a topos.)  Write
$\Open(X)$ for its poset of open subsets.  A presheaf on $X$ is a functor
$F\cln \Open(X)^\op \to \Set$.  It assigns to each open subset $U$ a set
$F(U)$, whose elements are called \demph{sections over $U$} (for reasons to be
explained).  It also assigns to each open $V \sub U$ a function $F(U) \to
F(V)$, called \demph{restriction} from $U$ to $V$ and denoted by $s \mapsto
s\restr{V}$.  I will write $\Pshf(X)$ for the category of presheaves on $X$.

\begin{examples}
\begin{enumerate}
\item   \label{eg:pshf-cts}
Let $F(U) = \{\text{continuous functions } U \to \reals\}$; restriction
is restriction.
\item   \label{eg:pshf-bdd}
The same, but with `bounded' in place of `continuous'.
\end{enumerate}
\end{examples}

Examples~\eqref{eg:pshf-cts} and~\eqref{eg:pshf-bdd} are qualitatively
different: continuity is a local property, but boundedness is not.  This
difference can be captured by asking the following question.  Let $(U_i)_{i
\in I}$ be a family of open subsets of $X$, and take, for each $i \in I$, a
section $s_i \in F(U_i)$.  Might there be some $s \in F(\bigcup_{i \in I}
U_i)$ such that $s\restr{U_i} = s_i$ for all $i$?

For this to stand a chance of being true, functoriality demands that the
sections $s_i$ must satisfy a `matching condition': $s_i\restr{U_i \cap U_j} =
s_j\restr{U_i \cap U_j}$ for all $i$ and $j$.  A \demph{sheaf} is a presheaf
such that for every family $(U_i)_{i \in I}$ of open sets and every matching
family $(s_i)_{i \in I}$, there is a unique $s \in F(\bigcup_{i \in I} U_i)$
such that $s\restr{U_i} = s_i$ for all $i \in I$.

\begin{examples}        \label{eg:sheaves}
\begin{enumerate}
\item The first example above, with continuous functions, is a sheaf.  The
proof can be split into two parts.  Given $(U_i)$ and $(s_i)$, there is
certainly a unique \emph{function} $s\cln \bigcup U_i \to \reals$ (continuous
or not) such that $s\restr{U_i} = s_i$ for all $i$.  The question now is
whether $s$ is continuous; and because continuity is a local property, it is.

\item The second example above, with bounded functions, is not a sheaf (for a
general space $X$).  This is because boundedness is \emph{not} a local
property.  

\item  \label{eg:shf-extend}
The sheaf of continuous real-valued functions is rather floppy, in the
sense that there are usually many ways to extend a continuous function from a
smaller set to a larger one.  Often people consider sheaves made up of
holomorphic or rational functions, which are much more rigid: there are
typically few or no ways to extend.  It is quite normal for there to be no
global sections (sections over $X$) at all.

\item   \label{eg:shf-bundle} 
Take any continuous map $Y \toby{p} X$ of topological spaces (which can be
thought of as a kind of bundle over $X$).  Then there arises a sheaf $F$ on
$X$, in which $F(U)$ is the set of continuous maps $s\cln U \to Y$ such that
the triangle on the left commutes:
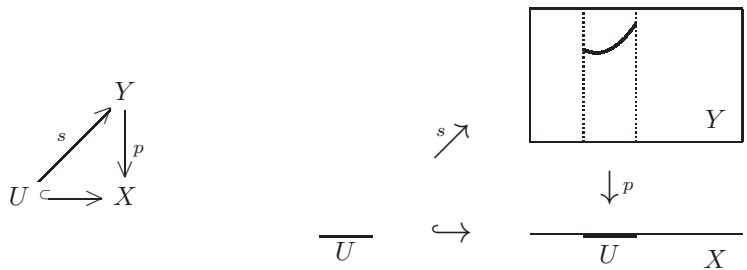
\[
\begin{diagram}
        &               &Y      \\
        &\ruTo<s        &\dTo>p \\
U       &\rIncl         &X      \\
\end{diagram}
\qquad\qquad\qquad
\setlength{\unitlength}{1em}
\setlength{\fboxsep}{0pt}
\begin{array}{c}
\begin{picture}(16,10)
\thicklines
\put(0,1.4){\line(1,0){2}}
\cell{1}{.5}{b}{U}
\put(10,1.4){\line(1,0){2}}
\cell{11}{.4}{b}{U}
\thinlines
\put(8,1.5){\line(1,0){8}}
\cell{15}{0.2}{b}{X}
\put(8,5){\framebox(8,5)}
\cell{15}{5.5}{b}{Y}
\qbezier[30](10,5)(10,7.5)(10,10)
\qbezier[30](12,5)(12,7.5)(12,10)
\qbezier(10,8.5)(11,8)(12,9.5)
\qbezier(10,8.45)(11,7.95)(12,9.45)
\cell{5}{1.5}{c}{\mbox{\Large\ensuremath{\incl}}}
\cell{11}{3.25}{c}{\mbox{\Large\ensuremath{\downarrow}}}
\cell{11.5}{3.25}{l}{\scriptstyle{p}}
\cell{5}{5}{c}{\mbox{\Large\ensuremath{\nearrow}}}
\cell{4.8}{5.2}{br}{\scriptstyle{s}}
\end{picture}%
\end{array}
\]
Such an $s$ is precisely a right inverse, or `section', of the map
$p^{-1}U \to U$ induced by $p$.
\end{enumerate}
\end{examples}

There is also an abstract categorical explanation of where the concept of
sheaf comes from.  Fix a space $X$.  We have a functor
\[
I\cln \Open(X) \to \TopSp/X
\]
where $\TopSp$ is the category of topological spaces, $\TopSp/X$ is the slice
category, and $I(U) = (U \incl X)$.  This functor $I$ embodies the simple
thought that an open subset of a topological space can be treated as a space
in its own right.  We now apply to $I$ two very general categorical
constructions, from which the sheaf concept will appear automatically.

First, purely because the domain of $I$ is small and the codomain has small
colimits, there is an induced adjunction
\[
\Pshf(X) = \Set^{\Open(X)^\op}
\pile{\rTo^{\dashbk \otimes I}\\ \dbot\\ \lTo_{\Hom(I, \dashbk)}}
\TopSp/X.
\]
The right adjoint is given by
\[
(\Hom(I, Y))(U)
=
\TopSp/X
\,
(I(U), Y)
\]
where $Y = \left(\vslob{Y}{p}{X}\right) \in \TopSp/X$ and $U \in \Open(X)$.
This is, in fact, the process described in
Example~\ref{eg:sheaves}(\ref{eg:shf-bundle}): the sheaf $F$ defined there is
$\Hom(I, Y)$.  The left adjoint can be described as a coend or colimit: for $F
\in \Pshf(X)$,
\[
F \otimes I     
= 
\int^U F(U) \times I(U) 
=
\Bigl(
\bigl(
\Colt{U, s}
U
\bigr)
\to
X
\Bigr)
\]
where the colimit is over all $U \in \Open(X)$ and $s \in F(U)$, and the
map from the colimit to $X$ is the canonical one.

Second, every adjunction restricts%
\label{p:adjn-eqv}
canonically to an equivalence between full subcategories: one consists of the
objects at which the unit of the adjunction is an isomorphism, and the other
of the objects at which the counit is an isomorphism.  Write the
equivalence obtained from the adjunction above as 
\[
\Sh(X) 
\pile{\rTo\\ \deqv\\ \lTo}
\Et(X).
\]
It can be shown that this $\Sh(X)$ is the same category of sheaves as before.
In this way, the notion of sheaf arises canonically from the very simple
functor $I\cln \Open(X) \to \TopSp/X$.  The notion of \'etale bundle also
arises canonically: \'etale bundles over $X$ are (by definition, if you like)
the objects of $\Et(X)$.  Among other things, this equivalence shows that
every sheaf is of the form described in
Example~\ref{eg:sheaves}(\ref{eg:shf-bundle}).  See \citet{MaMo} for details.

One way or another, we have the category $\Sh(X)$ of sheaves on $X$.  It is a
topos.  Its subobject classifier $\Omega$ is given by
\[
\Omega(U) 
=
\{
\textrm{open subsets of } U
\}.
\]
The crucial fact about $\Sh(X)$ is that---modulo a small lie that I will
repair later---
\slogan{$X$ can be recovered from $\Sh(X)$.}
So the class of topological spaces embeds into the class of toposes.  We can
think of toposes as generalized spaces.

A common technique in topos theory is to take a concept from topology or
geometry and extend it to toposes.  For example, suppose you hear someone
talking about `connected toposes'.  You may have no idea what one is, but you
can bet that the definition has been obtained by determining what
property of the topos $\Sh(X)$ corresponds to connectedness of the space
$X$, then taking that as the definition of connectedness for all toposes.

The next few subsections are all examples of this generalization process.

\chunk{Geometric morphisms}

So far I have said nothing about maps between toposes.  There is an obvious
candidate for what a map of toposes should be: a functor preserving finite
limits, exponentials, and subobject classifiers.  Such a functor is called a
\demph{logical morphism}.  They have a part to play, but there is another
notion of map of toposes that has been found much more useful.  It can be
derived by generalizing from topology.

Every map $f\cln X \to Y$ in $\TopSp$ induces an adjunction
\begin{equation}        \label{eq:Sh-GM}
\Sh(X)
\pile{\lTo^{f^*}\\ \dbot\\ \rTo_{f_*}}
\Sh(Y).
\end{equation}
This is not obvious.  The right adjoint $f_*$ is easy to construct---
\[
(f_* F)(V) = F(f^{-1} V)
\]
($F \in \Sh(X)$, $V \in \Open(Y)$)---but the left adjoint $f^*$ is harder.  It
can be made easy by invoking the equivalence between sheaves and \'etale
bundles; but I will not go into that, or give any other description of $f^*$.

It is a fact that $f^*$ preserves finite limits.  It is also a fact (modulo
the usual small lie) that there is a natural correspondence between continuous
maps $X \to Y$ and adjunctions~\eqref{eq:Sh-GM} in which the left adjoint
preserves finite limits.  So now we know what continuous maps look like in
topos-theoretic terms.  We duly generalize:

\begin{defn}    \label{defn:GM}
Let $\E$ and $\F$ be toposes.  A \demph{geometric morphism} $f\cln \E \to \F$
is an adjunction
\[
\E \pile{\lTo^{f^*}\\ \dbot\\ \rTo_{f_*}} \F
\]
in which the left adjoint $f^*$ preserves finite limits.  (People often say
`left exact left adjoint'.)  The right adjoint $f_*$ is called the
\demph{direct image} part of $f$, and $f^*$ is the \demph{inverse image} part.
\end{defn}

I will write $\Topos$ for the category of toposes and geometric morphisms.
(Really it's a 2-category, in an obvious way.)  By construction, we have a
functor
\[
\Sh\cln \TopSp \to \Topos
\]
which is (2-categorically) full and faithful, modulo the usual small lie.

\begin{examples}        \label{eg:GMs}
\begin{enumerate}
\item 
Every functor $f\cln \scat{C} \to \scat{D}$ induces a string of adjoint
functors
\[
\begin{diagram}[width=3em]
\Psh{\scat{C}} &
\pile{\rTo^{f_!}\\ \dbot\\ \lTo~{f^*}\\ \dbot\\ \rTo_{f_*}} &
\Psh{\scat{D}} \\
\end{diagram}
\]
between presheaf categories.  Here $f^* = \dashbk \of f$, and $f_!$ and $f_*$
are left and right Kan extension along $f$, respectively.  Since $f^*$ has a
left adjoint, it preserves limits.  Hence $(f^*, f_*)$ is a geometric morphism
$\Psh{\scat{C}} \to \Psh{\scat{D}}$.

\item   \label{eg:GM-sheaves}
It turns out that, for any topological space $X$, the inclusion $\Sh(X) \incl
\Pshf(X)$ has a finite-limit-preserving left adjoint.  It is called
\demph{sheafification} or the \demph{associated sheaf} functor.  So the
inclusion of sheaves into presheaves is a geometric morphism.

Since $\Sh(X)$ is a \emph{full} subcategory, the inclusion is full and
faithful; and for totally general reasons, this is equivalent to the counit of
the adjunction being an isomorphism.  In other words, sheafifying a sheaf does
not change it.
\end{enumerate}
\end{examples}

\chunk{Points}

Let us generalize another concept of topology.  The points of a topological
space $X$ correspond to the maps $1 \to X$ (where $1$ is the one-point space),
which correspond to the geometric morphisms $\Sh(1) \to \Sh(X)$.  But $\Sh(1)
= \Pshf(1) = \Set$, so we make the following definition.

\begin{defn}    \label{defn:point}
A \demph{point} of a topos $\E$ is a geometric morphism $\Set \to \E$.
\end{defn}

\chunk{Embeddings and Grothendieck toposes}

For any subspace $Y$ of a space $X$, the inclusion $Y \incl X$ is an
\demph{embedding}, that is, a homeomorphism to its image.  It can be shown
that a map $f\cln Y \to X$ of spaces is an embedding if and only if the direct
image part $f_*$ of the corresponding geometric morphism $f\cln \Sh(Y) \to
\Sh(X)$ is full and faithful.  So, as usual, we generalize:

\begin{defn}
A geometric morphism $f\cln \F \to \E$ is an \demph{embedding} (or
\demph{inclusion}) if the direct image functor $f_*$ is full and faithful.
\end{defn}

We then say that $\F$ is a \demph{subtopos} of $\E$.  At least, this is the
right thing to say up to equivalence.  Perhaps we should reserve that word for
when $\F$ is actually a (full) subcategory of $\E$ and $f_*$ is the inclusion
$\F \incl \E$, rather than allowing $f_*$ to be any old full and faithful
functor.  But a full and faithful functor induces an equivalence to its image,
so it makes no real difference.

Probably the easiest toposes are the \demph{presheaf toposes}: those
equivalent to $\Psh{\scat{C}} = \Set^{\scat{C}^\op}$ for some small category
$\scat{C}$.  So maybe subtoposes of presheaf toposes are relatively easy too.
They have a special name:

\begin{defn}    \label{defn:GT}
A topos is \demph{Grothendieck} if it is (equivalent to) a subtopos of some
presheaf topos.
\end{defn}

For instance, we saw in Example~\ref{eg:GMs}(\ref{eg:GM-sheaves}) that
$\Sh(X)$ is a subtopos of $\Pshf(X) = \Psh{\Open(X)}$, for any topological
space $X$.  Hence $\Sh(X)$ is a Grothendieck topos.

Being Grothendieck is generally thought of as a mild condition on a topos.  A
Grothendieck topos has all small limits, which immediately disqualifies
toposes such as $\FinSet$, $\FinSet^{\scat{C}^\op}$, etc.  But other than
toposes arising from finite sets (or sets subject to some other cardinality
bound), most of the toposes that people have worked with are Grothendieck.  A
notable exception is the effective topos, the maps in which can be thought of
as computable functions.  Other non-Grothendieck toposes occur in the
topos-theoretic approach to non-standard analysis.

There is a theorem of Giraud giving a list of conditions on a category
equivalent to it being a Grothendieck topos.  It includes non-elementary
axioms such as `there is a small generating set'.  (`Non-elementary' means
that it refers to a pre-existing notion of set.)  The Grothendieck toposes are
sometimes regarded as the nice toposes, but perhaps the definition of
Grothendieck topos is not as nice as the definition of elementary topos.

Definition~\ref{defn:GT} is not the definition of Grothendieck topos that you
will find in most books.  I will now give a brief indication of what the
standard definition is and why it is equivalent to the one above.

Fix a small category $\scat{C}$.  There is a one-to-one correspondence
between the subtoposes of $\Psh{\scat{C}}$ and the 
\demph{Grothendieck topologies} on $\scat{C}$.  A Grothendieck topology is a
kind of explicit, combinatorial structure; it specifies which diagrams
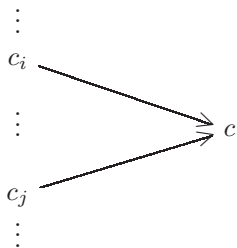
\[
\begin{diagram}[height=1.2em,width=4em]
\vdots  &               &       \\
c_i     &               &       \\
        &\rdTo(2,2)     &       \\
\vdots  &               &c      \\
        &\ruTo(2,2)     &       \\
c_j     &               &       \\
\vdots  &               &       \\
\end{diagram}
\]
in $\scat{C}$ are to be thought of as `covering families' and which are not.
(There are axioms.)  The motivating example is that given a topological space
$X$, there is a canonical Grothendieck topology on $\Open(X)$: a family $(U_i
\incl U)_{i \in I}$ of subsets of $U \in \Open(X)$ is covering if and only if
$U = \bigcup_{i \in I} U_i$.

The bijection
\[
\{ \textrm{Grothendieck topologies on } \scat{C} \}
\iso
\{ \textrm{subtoposes of } \Psh{\scat{C}} \}
\]
is written
\[
J \leftrightarrow \Sh(\scat{C}, J).
\]
A pair $(\scat{C}, J)$, consisting of a small category $\scat{C}$ equipped
with a Grothendieck topology $J$, is called a \demph{site}, and $\Sh(\scat{C},
J)$ is the category of \demph{sheaves} on that site.  For example, let $X$ be
a topological space, take $\scat{C} = \Open(X)$, and take $J$ to be the
Grothendieck topology mentioned above; then $\Sh(\scat{C}, J) = \Sh(X)$.  Most
books proceed as follows: define Grothendieck topology, define site, define
the category of sheaves on a site, then define a Grothendieck topos to be a
category equivalent to the category of sheaves on some site.

I do not know a short way to explain why the subtoposes of $\Psh{\scat{C}}$
correspond to the Grothendieck topologies on $\scat{C}$.  The following two
paragraphs may make it seem easier, or harder.

First, there is an explicit classification of the subtoposes of \emph{any}
topos $\E$.  Indeed, it can be shown that the subtoposes of $\E$ correspond to
the maps $j\cln \Omega \to \Omega$ satisfying certain equations.  
(Such a $j$ is called a \demph{Lawvere--Tierney topology} on $\E$, although
this is so distant from the original usage of the word `topology' that some
people object; Peter Johnstone, for instance, uses \demph{local operator}
instead.)  By definition of subobject classifier, it is equivalent to say that
a subtopos of $\E$ amounts to a subobject of $\Omega$ satisfying certain
axioms.

Second, take $\E = \Psh{\scat{C}}$.  We know
(Example~\ref{egs:toposes}(\ref{eg:topos-pshf})) that $\Omega \in
\Psh{\scat{C}}$ is given by $\Omega(c) = \{ \text{sieves on } c \}$.  Hence a
subtopos of $\Psh{\scat{C}}$ corresponds to a collection of sieves in
$\scat{C}$, satisfying certain axioms.  Calling these the `covering sieves'
gives the notion of Grothendieck topology.

\chunk{Locales}

Here I will explain the `small lie' mentioned several times above, and make
amends.  I will also explain why topos theorists are fond of jokes about
pointless topology.

The definition of sheaf on a topological space $X$ does not mention the points
of $X$.  It mentions only the open sets and inclusions between them, and uses
the fact that it is possible to take arbitrary unions and finite intersections
of open sets.  Having observed this, you can see why the space $X$ cannot
\emph{always} be recovered from the topos $\Sh(X)$.  For instance, if $X$ is
indiscrete (has no open sets except $\emptyset$ and $X$) and nonempty, then
$\Sh(X)$ is the same no matter how many points $X$ has.

The idea now is to split the process $X \goesto \Sh(X)$ into two steps.  First,
we forget the points of $X$, leaving just the set of open sets, ordered by
inclusion.  Then, we form the category of `sheaves' on that ordered set
(defined as for topological spaces, almost verbatim).

\begin{defn}
A \demph{frame} is a partially ordered set such that every subset has a join
($=$ least upper bound $=$ sup), every finite subset has a meet ($=$ greatest
lower bound $=$ inf), and finite meets distribute over joins.  A \demph{map of
frames} is a map preserving order, joins and finite meets.
\end{defn}

A topological space $X$ has a frame $\Open(X)$ of open subsets, and a
continuous map $f\cln X \to Y$ induces a map $f^{-1}\cln \Open(Y) \to
\Open(X)$ of frames.  This gives a functor
\[
\Open\cln \TopSp \to \Frame^\op.
\]
We now perform a linguistic manoeuvre.  $\Frame^\op$ is the desired category
of `pointless spaces'.  But we cannot wholeheartedly say that a frame is a
pointless space, because the \emph{maps} of frames are the wrong way round.
So we introduce a new word---\demph{locale}---and define the category $\Loc$
of locales by $\Loc = \Frame^\op$.  We can wholeheartedly say that a
\emph{locale} is a pointless space.

There is a functor $\Sh\cln \Loc \to \Topos$, defined just as for topological
spaces except that unions become joins and intersections become meets.  The
functor $\Sh\cln \TopSp \to \Topos$ factorizes as
\[
\TopSp \toby{\Open} \Loc \toby{\Sh} \Topos.
\]
This is the two-step process mentioned above.

Whenever I have said `modulo a small lie', you can interpret that as `use
locales instead of topological spaces'.  For example, $\Sh\cln \Loc \to
\Topos$ really is full and faithful, in a suitably up-to-isomorphism sense:
locale maps $X \to Y$ correspond one-to-one with isomorphism classes of
geometric morphisms $\Sh(X) \to \Sh(Y)$.  This means that $\Loc$ is equivalent
to a full subcategory of $\Topos$.  (Actually it is an equivalence of
2-categories, but I will gloss over that point.)

Every locale gives rise to a topos---but the converse is also true.  Given a
topos $\E$, the subobjects of $1$ form a poset $\Sub_\E(1)$.  Assuming that
$\E$ has enough colimits, $\Sub_\E(1)$ is a frame.  This process defines a
functor
\[
\begin{array}{ccc}
\Topos  &\to            &\Loc           \\
\E      &\goesto        &\Sub_\E(1).
\end{array}
\]
I am now quietly changing $\Topos$ to mean the toposes with small colimits;
this includes all Grothendieck toposes.

You might think that $1$ could have no interesting subobjects, since that is
the case in the most obvious topos, $\Set$.  But there are toposes that are
nearly as obvious in which $\Sub_\E(1)$ is not trivial.  For instance, take
$\E = \Set^I$ for any set $I$: then $\Sub_\E(1)$ is the power set of $I$.

Now a wonderful thing is true.  The functor just defined is left adjoint to
the inclusion $\Sh\cln \Loc \incl \Topos$.  This means that $\Loc$ is
(equivalent to) a \emph{reflective} subcategory of $\Topos$.  Hence the counit
is an isomorphism:
\[
X \iso \Sub_{\Sh(X)}(1)
\]
for any locale $X$.  This is how you recover a locale from its topos of
sheaves.

So $\Loc$ sits inside $\Topos$ as a subcategory of the best kind: full and
reflective, like abelian groups in groups.  It is reasonable to say
that a locale is a special sort of topos.  More formally, a topos is
\demph{localic} if it is of the form $\Sh(X)$ for some locale $X$.  Localic
toposes are easy to work with; if you were having trouble proving something
for arbitrary toposes, you might start by trying to prove it in this
special case.  

Since every locale is of the form $\Sub_\E(1)$ for some topos $\E$, locale
theory can be regarded as the fragment of topos theory concerning subobjects
of $1$.  A subobject of $1$ is a map $1 \to \Omega$, which can reasonably
called a truth value.  In that sense, locale theory is the study of truth
values. 

The notion of locale can also be seen as a decategorification of the notion of
Grothendieck topos.  A poset $P$ is a category enriched in the two-element
totally ordered set $2$.  There is a Yoneda embedding $P \to 2^{P^\op}$, which
has a finite-meet-preserving left adjoint if and only if $P$ is a frame.
Analogously, it is almost true that for a category $\E$, the Yoneda embedding
$\E \to \Set^{\E^\op}$ has a finite-limit-preserving left adjoint if and only
if $\E$ is a Grothendieck topos.  (This result is due to \citet{StrNT}.
`Almost' refers to a set-theoretic size condition.)  A map of frames is a
function preserving joins and finite meets, and the inverse image
part of a geometric morphism is a functor preserving colimits and finite
limits.  Thus, locales play roughly the same role among 2-enriched categories
as Grothendieck toposes play among $\Set$-enriched categories.

How much has been lost by passing from topological spaces to locales?  In most
people's view, not much.  For example, we observed that all nonempty
indiscrete spaces give rise to the same locale; but many mathematicians regard
indiscrete spaces with $\geq 2$ points as `pathological' and would be
positively happy to see them go.

In fact, some things are gained.  For example, a subgroup of a topological
group need not be closed, and non-closed subgroups are often regarded as
pathological (since the corresponding quotients are non-Hausdorff).  But it
is a theorem that every subgroup of a \emph{localic} group is closed.  See for
instance Section~C5.3 of \citet{JohSE}.

The functor $\Open\cln \TopSp \to \Loc$ has a right adjoint, which I will not
describe.  As mentioned on page~\pageref{p:adjn-eqv}, every adjunction
restricts canonically to an equivalence between full subcategories.  In this
case, this gives an equivalence between:
\begin{itemize}
\item a full subcategory of $\TopSp$, whose objects are called the
\demph{sober} spaces
\item a full subcategory of $\Loc$, whose objects are called the
\demph{spatial} locales.
\end{itemize}
Another way of interpreting the phrase `modulo a small lie' is `true for sober
spaces'.  Sobriety amounts to a rather mild separation condition.  For
example, every Hausdorff space is sober.  So in passing from a Hausdorff space
to a locale, or to a topos, nothing whatsoever is lost.

There is a kind of attitudinal paradox here.  Many algebraic topologists think
\emph{only} about Hausdorff spaces, and regard non-Hausdorff spaces as
pathological.  But these are often the same people who feel strongly that
topological spaces are not really about open sets; they think in terms of
points and paths and homotopies.  So it is perhaps paradoxical that the
Hausdorff condition guarantees that a space can be understood in terms of
its open sets alone: the topos of sheaves depends on nothing else, and
contains all the information about the original space.

\section{Toposes and universal algebra}
\label{sec:univ-alg}

The point of this section is to explain what people mean when they talk about
the classifying topos of a theory.  Another way to look at it is this: I will
explain how toposes can be viewed as cousins of operads and Lawvere theories.

In classical universal algebra, an algebraic theory (or strictly, a
presentation of an algebraic theory) consists of a bunch of operation symbols
of specified arities, together with a bunch of equations.  To take the
standard example, the (usual presentation of the) theory of groups consists of
\begin{itemize}
\item an operation symbol $1$ of arity $0$
\item an operation symbol $\blank^{-1}$ of arity $1$
\item an operation symbol $\cdot$ of arity $2$
\end{itemize}
together with the usual equations.  You can speak of `models' of an algebraic
theory in any category $\E$ with finite products.  In our example, they are
the internal groups in $\E$.

But there are other ways of looking at such theories.  

Consider the free finite product category $\T$ equipped with an internal
group.  (There are general reasons why such a thing must exist.)  Its
universal property is that for any finite product category $\E$, the
finite-product-preserving functors $\T \to \E$ correspond to the internal
groups in $\E$.

Concretely, $\T$ looks something like this.  It must contain an object $X$,
the underlying object of the internal group.  Since $\T$ has finite
products, it must also contain an object $X^n$ for each $n \in \nat$.  There
is no reason for it to have any other objects, and since it is free, it does
not.  A map $X^n \to X^m$ is (by definition of product) an $m$-tuple of maps
$X^n \to X$; and the maps $X^n \to X$ are (by freeness) whatever maps $G^n \to
G$ must exist for \emph{any} internal group $G$ in \emph{any} finite product
category.  That is, they are the $n$-ary operations in the theory of groups:
the words in $n$ letters.

This category $\T$ is called the \demph{Lawvere theory of groups}.  The same
goes for rings, lattices, etc.  In all these cases, $\T$ is a finite product
category with the further property that the objects are in bijection with the
natural numbers, the product of objects corresponding to addition of numbers.
This further property holds because the theories described so far have been
single-sorted: a model is a \emph{single} object equipped with some structure.

But there are also many-sorted theories, such as the two-sorted theory of
pairs $(R, M)$ in which $R$ is a ring and $M$ an $R$-module.  So we can widen
the notion of algebraic theory to include all (small) finite product
categories.  Some people say that an algebraic theory \emph{is} just a finite
product category.  Others say that algebraic theories \emph{correspond}
to finite product categories.  Others still, more traditionally, say that
algebraic theories correspond to only \emph{certain} finite product
categories.

Terminology aside, we can play the same game for other classes of limit.  For
example, it makes no sense to talk about internal categories in an arbitrary
finite product category, because the definition of internal category needs
pullbacks.  (Composition in an internal category $\scat{C}$ is a map
$\scat{C}_1 \times_{\scat{C}_0} \scat{C}_1 \to \scat{C}_1$.)  But we can talk
about internal categories in a finite \emph{limit} category; and as before,
there is a free finite limit category $\T$ equipped with an internal category.
This means that for any finite limit category $\E$, the
finite-limit-preserving functors $\T \to \E$ correspond to the internal
categories in $\E$.  A small finite limit category is called (or corresponds
to) an \demph{essentially algebraic theory}.

In a category with finite products you can talk about internal groups but not,
in general, internal categories.  In a category with finite limits you can
talk about both.  By extending the list of properties that the category is
assumed to satisfy, you can accommodate more and more sophisticated kinds of
theory.  (The theory of internal categories is more `sophisticated' than that
of groups in the sense that composition is only defined for \emph{some} pairs
of maps, whereas classical universal algebra can only handle operations
defined on \emph{all} pairs.)  The properties need not be of the form `limits
of such-and-such a type exist'.  For example, it is sometimes useful to assume
epi-mono factorization, as we shall see.

There is a trade-off here.  As you allow more sophisticated language, you
widen the class of theories that can be expressed, but you narrow the class of
categories in which it makes sense to take models.  (You also make more work
for yourself.)  In the same way, if you trade in your motorbike for a
double-decker bus, you increase the number of passengers you can carry, but
you restrict where you can carry them: no low bridges or tight alleyways.
(You also increase your fuel costs.)  It is sensible, then, to use
the smallest class of theories containing the ones you are interested in.
For example, you \emph{could} treat groups as an essentially algebraic theory,
but that would mean you could only take models in categories with \emph{all}
finite limits, when in fact just products would do.

Before I get onto toposes, I want to point out a slightly different direction
that you can take things in.  Rather than just altering the \emph{properties}
that the categories are assumed to have, you can also alter the
\emph{structure} with which they are equipped. 

Take monoidal categories, for instance.  We can speak of internal monoids in
any monoidal category.  Hence, the theory of monoids can be regarded as the
free monoidal category containing an internal monoid.  (This is in fact the
category of finite ordinals.)  Similarly, it makes sense to speak of algebras
for an operad $P$ in any monoidal category, and we can associate to $P$ the
free monoidal category $\T$ containing a $P$-algebra.  Thus, for any monoidal
category $\E$, monoidal functors $\T \to \E$ correspond to $P$-algebras in
$\E$.

We might define a \demph{monoidal theory} to be a small monoidal category.
This gets us into the territory of PROPs, where there are nontrivial theorems
such as the classification of 2-dimensional topological quantum field
theories: the symmetric monoidal theory of (or, `PROP for') commutative
Frobenius algebras is the category of smooth 1-manifolds and diffeomorphism
classes of cobordisms.

All of this is to give an impression of how far-reaching these ideas
are.  It is a sketch of the context in which classifying toposes can be
understood.

You will have guessed that the same kind of thing can be said for toposes as
for categories with finite products, finite limits, etc.  Since toposes have
very rich structure (\emph{much} more than just finite limits), they
correspond to a very wide class of theories indeed.  

An example of the kind of theory that can be interpreted in a topos is the
theory of fields.  (This is rather a feeble example, but I want to
keep it simple.)  A field is, of course, a commutative ring $R$ satisfying the
axioms 
\begin{equation}        \label{eq:field-nontriv}
0 \neq 1
\end{equation}
and
\begin{equation}        \label{eq:field-main}
\forall x \in R, 
\quad 
x = 0 \textrm{ or } \exists y\cln xy = 1.
\end{equation}
By a mechanical process, this definition can be turned into a definition of
`internal field in a topos'.  As compensation for the imprecision of the rest
of this section, I will give the definition in detail; but if you want to skip
it, the point to retain is that it \emph{is} a mechanical process.

Let $\E$ be a topos.  We certainly know how to define `commutative ring in
$\E$': that makes sense in any category with finite products.  Let $R$ be a
commutative ring in $\E$.  The nontriviality axiom, $0 \neq 1$, is expressed
by saying that 
the equalizer of
\[
1 \parpair{0}{1} R
\]
is the initial object $0$.  For the other axiom, let us first define the
subobject $U \monic R$ consisting of the units (invertible elements).  The
`set' $P = \{ (x, y) \such xy = 1 \}$ is the pullback
\[
\begin{diagram}
P\SEpbk         &\rTo           &1              \\
\dMono          &               &\dMono>1       \\
R \times R      &\rTo_\cdot     &R.
\end{diagram}
\]
Now we want to define the `set' $U$ of units as the image of the composite map 
\[
f 
= 
\left(
P \monic R \times R \toby{\pr_1} R
\right).
\]
We can talk about images in a topos, since every map in a topos
factorizes essentially uniquely as an epi followed by a mono.  So, define $U
\monic R$ by the factorization
\[
f 
=
\left(
P \epic U \monic R
\right).
\]
The second field axiom states that every element of $R$ lies in either the
subobject $1 \monicby{0} R$ or the subobject $U \monic R$.  In other words, it
states that the map
\[
1 + U \to R
\]
is epi.  Here we have used the fact that every topos has coproducts, written
$+$.  

If you have read Section~\ref{sec:sets}, you will recognize that the informal
talk of `sets' (really, objects of $\E$) and the use of set-theoretic notation
$\{ \ldots \such \ldots \}$ are something to do with the internal language of
a topos.  This gives a hint of how the process can be mechanized.

(There are actually \emph{several} possible theories of fields, depending on
exactly how you write down the axioms.  They all have the same models in
$\Set$---namely, fields---but they do not have the same models in other
toposes.  For example, a genuinely different theory is obtained by changing
axiom~\eqref{eq:field-main} to `$\forall x \in R$, $(\not\!\exists y: xy = 1)
\implies x = 0$'.  But this does not affect the main point: given a list of
formally-expressed axioms such as~\eqref{eq:field-nontriv}
and~\eqref{eq:field-main}, there is an automatic process converting it into a
definition that makes sense in an arbitrary topos.)

You now have the choice between a short story and a long story. 

The short story is that what we did for finite product and finite limit
categories can also be done for toposes.  The theories corresponding to
toposes are called the geometric theories, and the topos corresponding to a
particular geometric theory is called its classifying topos.

The long story is longer because there are two different notions of map of
toposes---and you need to decide what a map of toposes is in order to state
the universal property of the topos resulting from a theory.  

The more obvious but less used notion of map of toposes is a functor
preserving all the structure in sight: finite limits, exponentials, and the
subobject classifier.  These are called \demph{logical morphisms}.  Now in a
topos, you can interpret a really vast range of theories: any `higher-order
theory', in fact.  (First order means that you can only quantify over
\emph{elements} of a set; in a second order theory you can also quantify over
\emph{subsets} of a set; and so on.)  Models of any such theory get along well
with logical morphisms, because logical morphisms preserve everything.  So you
can tell a similar story for toposes, logical morphisms and higher order
theories as for finite product categories, finite-product-preserving functors
and algebraic theories.

The more popular notion of map of toposes is that of geometric
morphism.  (Here it helps to have read Section~\ref{sec:geom}, where the
definition is motivated.)  A \demph{geometric morphism} between toposes is a
functor with a finite-limit-preserving left adjoint.  The corresponding
theories are the \demph{geometric theories}.  I will not give the definition,
but it is not too bad an approximation to say that they are the same as the
first-order theories: every geometric theory is first-order, and almost every
first-order theory that one encounters is geometric.

Given a geometric theory, a \demph{classifying topos} for the theory is a
cocomplete topos $\T$ with the property that for any cocomplete topos $\E$,
models of the theory in $\E$ correspond naturally to geometric morphisms $\E
\to \T$.  Every geometric theory has a classifying topos.  

There are two surprises here.  One is the appearance of the word `cocomplete',
which I will not explain and will not bother inserting below.  It is generally
thought of as a mild condition (satisfied by any Grothendieck topos, for
instance).  

The bigger surprise is the reversal of direction.  The previous cases lead us
to expect models in $\E$ to correspond to maps $\T \to \E$.  However, since a
geometric morphism is a pair of adjoint functors, the choice of direction is a
matter of convention.  As the name suggests, the choice that society made was
motivated by geometry.  Perhaps if the motivation had been universal algebra,
it would have been the other way round.  (This is an aspect of the thought
that geometry is dual to algebra.)  A map of toposes would then have been a
finite-limit-preserving functor with a right adjoint, which is more or less
the same thing as a functor preserving finite limits and small colimits.

If a topos is thought of as a generalized space (as in Section~\ref{sec:geom})
then the classifying topos of a theory can be thought of as its space of
models.  Indeed, a point of the classifying topos $\T$ is (by
Definition~\ref{defn:point}) a geometric morphism $\Set \to \T$, which is
exactly a model of the theory in $\Set$.  Some familiar topological spaces can
be construed as classifying toposes.  For example, there is a `theory of
Dedekind cuts' whose classifying topos is $\Sh(\reals)$, that is, $\reals$
regarded as a topos.

Given how much structure a topos contains, it is surprising how many
classifying toposes can be described simply.  I will now describe the
classifying topos of any algebraic theory, by the venerable expository device
of doing it just for groups.

We will need the notion of finite presentability.  A group (in $\Set$) is
\demph{finitely presentable} if it admits a presentation by a finite set of
generators subject to a finite set of relations.  The category of finitely
presentable groups and all homomorphisms between them will be written
$\Gp_\fp$.

\paragraph*{Aside} Finite presentability is a more categorical concept than it
might seem.  Writing $T\cln \Set \to \Set$ for the free group monad, a relation
(equation) in a set $X$ of generators is an element of $TX \times TX$.  So, a
family $(r_i)_{i \in I}$ of relations is a map $I \to TX \times TX$, or
equivalently a diagram 
\[
I \parpairu TX
\]
in $\Set$, or equivalently a diagram 
\[
FI \parpairu FX
\]
in $\Gp$, where $F\cln \Set \to \Gp$ is the free group functor.  The group
presented by these generators and relations is the coequalizer of this diagram
in $\Gp$.  Hence a group is finitely presentable precisely when it is the
coequalizer of some diagram $FI \parpairu FX$ in which $I$ and $X$ are finite
sets.  

This formulation of finite presentability in $\Gp$ uses the free group functor
$F$.  But in fact, there is a general definition of finite presentability of
an object of any category.  I will not go into this.

\bigskip

As promised, the classifying topos for groups is easy to describe:
\begin{thm}
The classifying topos for groups is $\Set^{\Gp_\fp}$.
\end{thm}
In other words, for any topos $\E$, a group in $\E$ is the same thing as a
geometric morphism $\E \to \Set^{\Gp_\fp}$.  

The same goes for other algebraic theories.  This yields something
interesting even for very trivial theories.  Take the theory of objects, whose
models in a category $\E$ are simply objects of $\E$.  A finitely presentable
set is just a finite set.  Hence for any topos $\E$, objects of $\E$
correspond to geometric morphisms $\E \to \Set^{\FinSet}$.  The topos
$\Set^\FinSet$ is therefore called the \demph{object classifier}.

We have been asking, for a given theory, `what topos classifies it?'  But we
can turn the question round and ask, for a given topos $\T$, `what does $\T$
classify?'  In other words, what are the geometric morphisms from an arbitrary
topos $\E$ into $\T$?  It is a fact that \emph{every} topos $\T$ is the
classifying topos of some geometric theory---although given how wide a class
of theories that is, perhaps this does not say very much.

There are clean answers to this reversed question for many toposes $\T$.  In
particular, this is so when $\T$ is the topos $\Sh(\scat{C}, J)$ of sheaves on
a site (Section~\ref{sec:geom}).  Here I will just tell you the answer for a
smaller class of toposes.

\begin{thm}
Let $\scat{C}$ be a category with finite limits.  Then the presheaf topos
$\Psh{\scat{C}}$ classifies finite-limit-preserving functors out of
$\scat{C}$. 
\end{thm}
In other words, for any topos $\E$, a geometric morphism $\E \to
\Psh{\scat{C}}$ is the same thing as a finite-limit-preserving functor
$\scat{C} \to \E$.

(If you know about flat functors, you can drop the assumption that $\scat{C}$
has finite limits: for any small category $\scat{C}$, the presheaf topos
$\Psh{\scat{C}}$ classifies flat functors out of $\scat{C}$.  This is
one version of Diaconescu's Theorem.)

So there is a back-and-forth translation between geometric theories and the
toposes that classify them.  In many cases, this translation is surprisingly
straightforward.


\begin{thebibliography}{7}

\bibitem[Johnstone(2003)]{JohSE}
Johnstone, P.~T., 2003.
\newblock Sketches of an Elephant: A Topos Theory Compendium.
\newblock Oxford Logic Guides. Oxford University Press.

\bibitem[Lawvere(1964)]{LawETCS}
Lawvere, F.~W., 1964.
\newblock An elementary theory of the category of sets.
\newblock \emph{Proceedings of the National Academy of Sciences of the U.S.A.}
  \textbf{52}:1506--1511.
\newblock Reprinted as \emph{Reprints in Theory and Applications of Categories}
  12:1--35, 2005.

\bibitem[Lawvere and Rosebrugh(2003)]{LaRo}
Lawvere, F.~W. and R.~Rosebrugh, 2003.
\newblock Sets for Mathematics.
\newblock Cambridge University Press, Cambridge.

\bibitem[Mac~Lane and Moerdijk(1994)]{MaMo}
Mac~Lane, S. and I.~Moerdijk, 1994.
\newblock Sheaves in Geometry and Logic.
\newblock Springer, Berlin.

\bibitem[Mc{L}arty(2004)]{McLECS}
Mc{L}arty, C., 2004.
\newblock Exploring categorical structuralism.
\newblock \emph{Philosophia Mathematica} \textbf{12}:37--53.

\bibitem[Par{\'e}(1974)]{Pare}
Par{\'e}, R., 1974.
\newblock Colimits in topoi.
\newblock \emph{Bulletin of the American Mathematical Society}
  \textbf{80}:556--561.

\bibitem[Street(1981)]{StrNT}
Street, R., 1981.
\newblock Notions of topos.
\newblock \emph{Bulletin of the Australian Mathematical Society}
  \textbf{23}:199--208.

\end{thebibliography}
\end{document}